\input amstex\documentstyle{amsppt}  
\pagewidth{12.5cm}\pageheight{19cm}\magnification\magstep1
\topmatter
\title Generic character sheaves on groups over $\bold k[\epsilon]/(\epsilon^r)$ \endtitle
\author G. Lusztig\endauthor
\address{Department of Mathematics, M.I.T., Cambridge, MA 02139}\endaddress
\thanks{Supported in part by National Science Foundation grant 1303060.}\endthanks   
\endtopmatter   
\document

\define\tcx{\ti{\cx}}

\define\hZ{\hat Z}
\define\hK{\hat K}

\define\si{\sim}

\define\qua{\quad}

\define\op{\oplus}
   
\define\part{\partial}
\define\emp{\emptyset}

\define\ra{\rangle}
\define\n{\notin}

\define\m{\mapsto}
\define\do{\dots}
\define\la{\langle}
\define\bsl{\backslash}

\define\lra{\leftrightarrow}

\define\sub{\subset}    
\define\bxt{\boxtimes}
\define\T{\times}
\define\ti{\tilde}
\define\nl{\newline}
\redefine\i{^{-1}}
\define\fra{\frac}
\define\un{\underline}

\define\ot{\otimes}
\define\bbq{\bar{\QQ}_l}

\define\Ad{\text{\rm Ad}}
\define\Hom{\text{\rm Hom}}

\define\ind{\text{\rm ind}}

\define\tr{\text{\rm tr}}

\define\a{\alpha}
\redefine\b{\beta}

\redefine\d{\delta}
\define\e{\epsilon}

\define\io{\iota}

\define\p{\pi}
\define\ph{\phi}
\define\ps{\psi}
\define\r{\rho}
\define\s{\sigma}
\redefine\t{\tau}

\redefine\l{\lambda}
\define\z{\zeta}
\define\x{\xi}

\redefine\D{\Delta}

\define\Si{\Sigma}

\define\kk{\bold k}

\define\FF{\bold F}

\define\QQ{\bold Q}

\define\ZZ{\bold Z}

\define\cc{\Cal C}
\define\cd{\Cal D}
\define\ce{\Cal E}
\define\cf{\Cal F}
\define\cg{\Cal G}
\define\ch{\Cal H}

\define\cl{\Cal L}

\define\co{\Cal O}

\define\car{\Cal R}

\define\ct{\Cal T}

\define\cz{\Cal Z}
\define\cx{\Cal X}
\define\cy{\Cal Y}

\define\fb{\frak b}

\define\fg{\frak g}

\define\fn{\frak n}

\define\ft{\frak t}

\define\tG{\ti G}

\define\tX{\ti X}

\define\sha{\sharp}

\define\tce{\ti\ce}

\define\DL{DL}
\define\GE{Ge}
\define\LA{La}
\define\GRE{L1}
\define\CS{L2}
\define\GEN{L3}

\head Introduction\endhead
\subhead 0.1\endsubhead
Let $\kk$ be an algebraic closure of the finite field $\FF_q$ with $q$ elements where $q$ is a power of a
prime number $p$. Let $G$ be a connected reductive group over $\kk$ with a fixed split $\FF_q$-rational
structure, a fixed Borel subgroup $B$ defined over $\FF_q$, with unipotent radical $U$ and a fixed maximal 
torus $T$ of $B$ defined over $\FF_q$. Let $\fg,\fb,\ft,\fn$ be the Lie algebras of $G,B,T,U$.
We fix a prime number $l\ne p$. If $\l:T(\FF_q)@>>>\bbq^*$ is a 
character, we can lift $\l$ to a character $\ti\l:B(\FF_q)@>>>\bbq^*$ trivial on $U(\FF_q)$ and we can form 
the induced representation $\ind_{B(\FF_q)}^{G(\FF_q)}\ti\l$ of $G(\FF_q)$. Its character is the class 
function $G(\FF_q)@>>>\bbq$ given by
$$y\m\sum\Sb B(\FF_q)x\in B(\FF_q)\bsl G(\FF_q);\\ xyx\i\in B(\FF_q)\endSb\ti\l(xyx\i).\tag a$$
This class function has a geometric analogue. Namely, we consider the diagram $T@<h<<\tG@>\p>>G$
where $\tG=\{(Bx,y)\in(B\bsl G)\T G;xyx\i\in B\}$, $\p(Bx,y)=y$ is the Springer map and $h(Bx,y)=d(xyx\i)$;
here $d:B@>>>T$ is the obvious homomorphism with kernel $U$. 
Let $\ce$ be a fixed $\bbq$-local system of rank $1$ on $T$ such that $\ce^{\ot m}\cong\bbq$ for some 
$m\ge1$ prime to $p$.
The geometric analogue of (a) is the complex $L_1=\p_!h^*\ce\in\cd(G)$. (For any algebraic
variety $X$ over $\kk$, $\cd(X)$ denotes the bounded derived category of constructible $\bbq$-sheaves on 
$X$.) 
When $\ce$ is defined over $\FF_q$ and has characteristic function $\l$ then $L_1$ is defined over 
$\FF_q$ and its characteristic function is (up to a nonzero scalar factor) the function (a). Thus $L_1$ can 
be viewed as a categoryfied version of the function (a). More precisely, $L_1$ is (up to shift) a perverse
sheaf on $G$; indeed, one of the main observations of \cite{\GRE} was that the (proper) map $\p$ is small, 
which implies that $L_1$ is an intersection cohomology complex; this was the starting point of the theory of
character sheaves on $G$, see \cite{\CS}. This point of view is useful since the complexes $L_1$ are defined
independently of the $\FF_q$-structure and from them one can extract not only the characters (a) for any $q$
but even their twisted versions defined in \cite{\DL}. 

\subhead 0.2\endsubhead
For any integer $r\ge1$ we consider the ring $\kk_r=\kk[\e]/(\e^r)$ ($\e$ is an indeterminate). Let 
$G_r=G(\kk_r)$ be the group of points of $G$ with values in $\kk_r$, viewed as an algebraic group over 
$\kk$ of dimension $r\D$ where $\D=\dim G$. Let $B_r=B(\kk_r)$, $T_r=T(\kk_r)$, $U_r=U(\kk_r)$. Note that 
$G_r$ inherits from $G$ a natural $\FF_q$-structure and that $B_r,T_r,U_r$ are defined over $\FF_q$.
For $r=1$, $G_r$ reduces to $G$; we would like to extend as much as possible the results in 0.1 from $r=1$
to a general $r$. If $\l:T_r(\FF_q)@>>>\bbq^*$ is a character, we can lift $\l$ to a character 
$\ti\l:B_r(\FF_q)@>>>\bbq^*$ trivial on $U_r(\FF_q)$ and we can form 
the induced representation $\ind_{B(\FF_q)}^{G(\FF_q)}\ti\l$ of $G(\FF_q)$. Its character is the class 
function $G_r(\FF_q)@>>>\bbq$ given by
$$g'\m\sum\Sb B_r(\FF_q)g\in B_r(\FF_q)\bsl G_r(\FF_q);\\ gg'g\i\in B_r(\FF_q)\endSb\ti\l(gg'g\i).\tag a$$
It generalizes the function 0.1(a). Again this class function has a geometric analogue. Namely, we consider 
the diagram $T_r@<h_r<<\tG_r@>\p_r>>G_r$ where 
$$\tG_r=\{(B_rg,g')\in(B_r\bsl G_r)\T G_r;gg'g\i\in B_r\},$$
$$\p_r(B_rg,g')=g',h_r(B_rg,g')=d_r(gg'g\i);$$
here $d_r:B_r@>>>T_r$ is the obvious homomorphism with kernel $U_r$. We can identify $T=\cy\ot\kk^*$ where
$\cy$ is the lattice of one parameter subgroups of $T$ and $T(\kk_r)=\cy\ot\kk_r^*$ where $\kk_r^*$ is the
group of units of $\kk_r$. The isomorphism $\kk^*\T\kk^{r-1}@>\si>>\kk_r^*$, 
$$(a_0,a_1,\do,a_{r-1})\m a_0+a_1\e+\do+a_{r-1}\e^{r-1},$$
identifies $T(\kk_r)$ with $T\T\ft^{r-1}$.
Let $f_1,\do,f_{r-1}$ be linear functions $\ft@>>>\kk$ and let $\ce$ be as in 0.1. 
We can form the local system
$\ce\bxt\cl_{f_1}\bxt\do\bxt\cl_{f_{r-1}}$ on $T\T\ft^{r-1}=T(\kk_r)$ (for the notation $\cl_{f_i}$ see 0.3).
The geometric analogue of (a) is the complex 
$L=\p_{r!}h_r^*(\ce\bxt\cl_{f_1}\bxt\do\bxt\cl_{f_{r-1}})\in\cd(G_r)$.
Again from the complexes $L$ one can extract the characters (a) for any $q$. 
In this paper we are interested in the conjecture in \cite{\GEN, 8(a)} according to which, when $r\ge2$, $L$
is (up to shift) 
an intersection cohomology complex on $G_r$, provided that $f_{r-1}$ is sufficiently general. This
would imply that there is a theory of generic character sheaves on $G_r$.
The conjecture was proved in \cite{\GEN, no.12} in the case where $G=GL_2$ and $r=2$.

In this paper we give a method to attack the conjecture for any $G$ and even $r$ (but with some restriction
on $p$); we carry out the method in detail in the cases where $r=2$ and $r=4$ and we prove the conjecture 
in these cases (with some restriction on $p$). We also prove a weak form of the conjecture assuming that 
$r=3$ (see Theorem 4.7). I believe that the method of this paper should be applicable with any $r\ge2$.

Our method is to first replace $L$ by another complex $K$ which is a geometric (categorified) form of
the character of a representation constructed by G\'erardin \cite{\GE} in 1975, then to try to describe 
explicitly the Fourier-Deligne transform of $K$ on $G_r$ (viewed as a vector bundle over $G$). For $r=2$ and
$r=4$ we show that this is a simple perverse sheaf of a very special kind, namely one associated to a local 
system of rank $1$ on a closed smooth irreducible subvariety of $G_r$; by a result of Laumon, this implies
that $K$ is itself a simple perverse sheaf, up to twist. Finally, we show that $L$ is a shift 
of $K$ for the values of $r$ that we consider and this gives the desired result.

I wish to thank Dongkwan Kim for pointing out an inaccuracy in an earlier version of this paper.

\subhead 0.3. Notation\endsubhead
In this paper all algebraic varieties are over $\kk$. 
We fix a nontrivial homomorphism $\ps:\FF_p@>>>\bbq^*$. For any morphism $f:X@>>>\kk$ let 
$X_f=\{(x,\l)\in X\T\kk;\l^q-\l=f(x)\}$ and let $\io:X_f@>>>X$ be the Artin-Schreier covering 
$(x,\l)\m x$. Then $\io_!\bbq$ is a local system with a natural action of $\FF_p$ (coming from the 
$\FF_p$-action $\z:(x,\l)\m(x,\l+\z)$ on $\tX$); we denote by $\cl_f$ the $\ps$-eigenspace of this action 
(a local system of rank $1$ on $X$).

Let $\d=\dim T$. 

For $x\in G$ if $X$ is an element of $\fg$ or a subset of $\fg$ we write ${}^xX$ instead of $\Ad(x)X$ and 
${}_xX$ instead of $\Ad(x\i)X$. 

\head 1. The complex $K$\endhead
\subhead 1.1\endsubhead
Let $X_1,X_2,\do$ and $Y_1,Y_2,\do$ be two sequences of noncommuting indeterminates. From the 
Campbell-Baker-Hausdorff formula we deduce the equality 
$$(e^{\e X_1}e^{\e^2X_2}\do)(e^{\e Y_1}e^{\e^2Y_2}\do)=e^{\e z_1}e^{\e^2z_2}\do$$
where $z_i=z_i(X_1,X_2,\do,X_i,Y_1,Y_2,\do,Y_i)$, ($i\ge1$) are universal Lie polynomials with
coefficients in $\ZZ[(i!)\i]$. (Here $\e$ commutes with each $X_i,Y_i$.) For example,

$z_1(X_1,Y_1)=X_1+Y_1$,

$z_2(X_1,X_2,Y_1,Y_2)=X_2+Y_2+[X_1,Y_1]/2$,

$z_3(X_1,X_2,X_3,Y_1,Y_2,Y_3)=X_3+Y_3+[X_2,Y_1]-[X_1,[X_1,Y_1]]/6-[Y_1,[X_1,Y_1]]/3$.
\nl
We deduce that if $X_1,X_2,\do$, $X'_1,X'_2,\do$ and $Y_1,Y_2,\do$ are three sequences of noncommuting 
indeterminates then we have the equality
$$(e^{\e X'_1}e^{\e^2X'_2}\do)(e^{\e Y_1}e^{\e^2Y_2}\do)(e^{\e X_1}e^{\e^2X_2}\do)\i
=e^{\e u_1}e^{\e^2u_2}\do$$
where $u_i=u_i(X'_1,\do,X'_i,Y_1,\do,Y_i,X_1,\do,X_i)$, ($i\ge1$) are universal Lie polynomials
with coefficients in $\ZZ[(i!)\i]$. For example,
$$u_1(X_1,Y_1,X'_1)=X'_1-X_1+Y_1,$$
$$u_2(X_1,X_2,Y_1,Y_2,X'_1,X'_2)=X'_2-X_2+Y_2+[X'_1,Y_1]/2-[X'_1,X_1]/2-[Y_1,X_1]/2,$$
$$\align&u_3(X_1,X_2,X_3,Y_1,Y_2,Y_3,X'_1,X'_2,X'_3)\\&=
X'_3-X_3+Y_3+[X'_2,Y_1]+[X_2,X_1]-[X'_2,X_1]-[Y_2,X_1]-[X'_1,[X'_1,Y_1]]/6\\&
-[Y_1,[X'_1,Y_1]]/3+[X_1,[X'_1,Y_1]]/2+[X'_1,[X'_1,X_1]]/6+[X'_1,[Y_1,X_1]]/6 \\&
+[Y_1,[X'_1,X_1]]/6+[Y_1,[Y_1,X_1]]/6-[X_1,[X'_1,X_1]]/3-[X_1,[Y_1,X_1]]/3.\endalign$$
Note that

(a) $u_i(X'_1,\do,X'_i,Y_1,\do,Y_i,X_1,\do,X_i)=X'_i-X_i+Y_i+u'_i$ where 

$u'_i=u'_i(X'_1,\do,X'_{i-1},Y_1,\do,Y_{i-1},X_1,\do,X_{i-1})$
\nl
is a Lie polynomial in $X'_1,\do,X'_{i-1},Y_1,\do,Y_{i-1},X_1,\do,X_{i-1}$.

\subhead 1.2\endsubhead
We now fix $r\ge2$. We write $r=2r'$ if $r$ is even and $r=2r'+1$ if $r$ is odd. We always assume that 
$p\ge r$. Then for any $X\in\fg$ and any $m\ge1$, the exponential 
$e^{\e^m X}\in G_r$ is well defined. For any $X_1,X_2,\do,X_{r-1}$ in $\fg$ we set
$$|X_1,X_2,\do,X_{r-1}|=e^{\e X_1}e^{\e^2X_2}\do e^{\e^{r-1}X_{r-1}}\in G_r.$$
We have an isomorphism of algebraic varieties 
$$G\T\fg^{r-1}@>\si>>G_r$$
given by 
$$(x,X_1,X_2,\do,X_{r-1})\m x|X_1,X_2,\do,X_{r-1}|=|{}_xX_1,{}_xX_2,\do,{}_xX_{r-1}|x.$$  
This restricts to isomorphisms of algebraic varieties $B\T\fb^{r-1}@>\si>>B_r$, $U\T\fn^{r-1}@>\si>>U_r$, 
$T\T\ft^{r-1}@>\si>>T_r$. (The last isomorphism is the same as one in 0.2.)

Let $X_1,X_2,\do,X_{r-1}$ and $Y_1,Y_2,\do,Y_{r-1}$ be two sequences in $\fg$ and let $x,y$ be in $G$. We 
have
$$(x|X_1,\do,X_{r-1}|)(y|Y_1,\do,Y_{r-1}|)=xy|Z_1,\do,Z_{r-1}|$$
where $Z_i=z_i({}_yX_1,\do,{}_yX_i,Y_1,\do,Y_i)\in\fg$ ($i=1,\do,r-1$) with notation of 1.1 and 
where $[,]$ becomes the Lie bracket in $\fg$; note that $Z_i$ are well defined since $p\ge r$. Moreover, we
have
$$(x|X_1,\do,X_{r-1}|)(y|Y_1,\do,Y_{r-1}|)(x|X_1,\do,X_{r-1}|)\i=xyx\i|U_1,\do,U_{r-1}|$$
where $U_i={}^xu_i({}_yX_1,\do,{}_yX_i,Y_1,\do,Y_i,X_1,\do,X_i)\in\fg$ ($i=1,\do,r-1$) with notation of 
1.1); note that $U_i$ are well defined since $p\ge r$. 

\subhead 1.3\endsubhead
Let $\ph:E@>>>X$ be an algebraic vector bundle with fibres of constant dimension $N$. Let $f:E@>>>\kk$ be a
morphism such that for any $x\in X$ the restriction $f^x:\ph\i(x)@>>>\kk$ is affine linear. Let $X_0$ be the
set of all $x\in X$ such that $f^x$ is a constant (depending of $x$) and let $f_0:X_0@>>>\kk$ be such that
$f(e)=f_0(\ph(e))$ for all $e\in\ph\i(X_0)$. Let $j:X_0@>>>X$ be the (closed) imbedding. We show:

(a) $\ph_!\cl_f\cong j_!\cl_{f_0}[-2N]$.
\nl
For any $x\in X-X_0$ we have $H^i_c(\ph\i(x),\cl_f)=0$ for all $i$. Hence $\ph_!\cl_f|_{X-X_0}=0$. We are 
reduced to the case where $X=X_0$. In this case we have $\cl_f=\ph^*\cl_{f_0}$ hence
$$\ph_!\cl_f=\ph_!\ph^*\cl_{f_0}=\cl_{f_0}\ot\ph_!\ph^*\bbq\cong\cl_{f_0}[-2N],$$
as required. (We ignore Tate twists.)

If in addition we are given a local system $\cf$ on $X$ and we denote $\ph^*\cf$ and $j^*\cf$ again by 
$\cf$, then from (a) we have immediately

(b) $\ph_!(\cf\ot\cl_f)\cong j_!(\cf\ot\cl_{f_0})[-2N]$.

\subhead 1.4\endsubhead
In the rest of this paper we assume that a nondegenerate symmetric bilinear invariant form 
$\la,\ra:\fg\T\fg@>>>\kk$ is given and that a sequence $A_1,A_2,\do,A_{r-1}$ of elements of $\ft$ is given 
such that $A_{r-1}$ is regular semisimple. This requires a further restriction on $p$ in addition to the 
restriction $p\ge r$. 

For a subspace $E$ of $\fg$ we set $E^\perp=\{\x\in\fg;\la \x,E\ra=0\}$.

Let $\cx$ be the variety of all
$$(Tx,y,X_1,X_2,\do,X_{r-1},Y_1,Y_2,\do,Y_{r-1})\in(T\bsl G)\T G\T\fg^{2r-2}$$
such that $xyx\i\in T$ and
$$u_j({}_yX_1,\do,{}_yX_j,Y_1,\do,Y_j,X_1,\do,X_j)\in {}_x\ft\text{ for }1\le j\le r'-1,$$
$$u_j({}_yX_1,\do,{}_yX_j,Y_1,\do,Y_j,X_1,\do,X_j)\in {}_x\fb\text{ if $j=r'$ and $r$ is odd}.$$
We have a diagram
$$G_r@<\p<<\cx@>h>>\kk$$
where $\p(Tx,y,X_1,\do,X_{r-1},Y_1,\do,Y_{r-1})=y|Y_1,\do,Y_{r-1}|$,
$$\align&h(Tx,y,X_1,\do,X_{r-1},Y_1,\do,Y_{r-1})\\&
=\sum_{j\in[1,r-1]}\la {}_xA_j,u_j({}_yX_1,\do,{}_yX_j,Y_1,\do,Y_j,X_1,\do,X_j)\ra.\endalign$$
(Note that if $x$ is replaced by $tx$, ($t\in T$) in the last sum, the sum remains unchanged since 
${}_tA_j=A_j$ for all $j$.) 
We define $\io:\cx@>>>T$ by
$$\io(Tx,y,X_1,X_2,\do,X_{r-1},Y_1,Y_2,\do,Y_{r-1})=xyx\i$$
and we set $\tce=\io^*\ce$. Let $K=\p_!(\tce\ot\cl_h)\in\cd(G_r)$.
Via the identification $G_r=G\T\fg^{r-1}$ (see 1.2) we can regard $G_r$ as a vector bundle over
$G$ with fibre $\fg^{r-1}$ endowed with a nondegenerate symmetric bilinear form. Hence the Fourier-Deligne
transform $\hK\in\cd(G_r)$ along these fibres is well defined. More explicitly, for $i=1,2$ we have the
 diagram $G_r@<\r_i<<G_r\T_G G_r@>h'>>\kk$ where $\r_i$ is the projection to the $i$-th factor and
$$h'(x|Y_1,\do,Y_{r-1}|,x|R_1,R_2,\do,R_{r-1}|)=\sum_{j\in[1,r-1]}\la Y_j,R_j\ra.$$
Then $\hK=\r_{2!}(\r_1^*K\ot\cl_{h'})[(r-1)\D]$ that is,
$$\hK=\r_{2!}(\r_1^*\p_!(\tce\ot\cl_h)\ot\cl_{h'})[(r-1)\D].$$
Let $\tcx$ be the variety of all
$$(Tx,y,X_1,X_2,\do,X_{r-1},Y_1,Y_2,\do,Y_{r-1},R_1,R_2,\do,R_{r-1})\in(T\bsl G)\T G\T\fg^{3r-3}$$
such that $xyx\i\in T$ and
$$u_j({}_yX_1,\do,{}_yX_j,Y_1,\do,Y_j,X_1,\do,X_j)\in{}_x\ft\text{ for }1\le j\le r'-1,$$
$$u_j({}_yX_1,\do,{}_yX_j,Y_1,\do,Y_j,X_1,\do,X_j)\in{}_x\fb\text{ if $j=r'$ and $r$ is odd }.$$
We have a cartesian diagram
$$\CD
\tcx@>\ti\r_1>>\cx\\
@V\s VV      @V\p VV\\
G_r\T_G G_r@>\r_1>>G_r
\endCD$$
where 
$$\align&\ti\r_1(Tx,y,X_1,\do,X_{r-1},Y_1,\do,Y_{r-1},R_1,\do,R_{r-1})\\&=
(Tx,y,X_1,\do,X_{r-1},Y_1,\do,Y_{r-1}),\endalign$$
$$\align&\s(Tx,y,X_1,\do,X_{r-1},Y_1,\do,Y_{r-1},R_1,\do,R_{r-1})\\&
=(y|Y_1,\do,Y_{r-1}|,y|R_1,\do,R_{r-1}|).\endalign$$
It follows that
$$\align&\hK=\r_{2!}(\s_!(\tce\ot\ti\r_1^*\cl_h\ot\cl_{h'})[(r-1)\D]\\&=   
\r_{2!}\s_!(\tce\ot\cl_{\ti h'}\ot\cl_{\ti h''})=(\r_2\s)_!(\tce\ot\cl_{\ti h})[(r-1)\D]\endalign$$
where $\ti h''=h\tr_1:\tcx@>>>\kk$, $\ti h'=h'\s:\tcx@>>>\kk$, $\ti h=\ti h'+\ti h'':\tcx@>>>\kk$ and the
inverse image of $\tce$ under $\tcx@>>>T$,
$$(Tx,y,X_1,\do,X_{r-1},Y_1,\do,Y_{r-1},R_1,\do,R_{r-1})\m xyx\i$$
is denoted again by $\tce$. Thus, 
$$\hK=\p'_!(\tce\ot\cl_{\ti h})[(r-1)\D]$$ 
where $\p':\tcx@>>>G_r$ and $\ti h:\tcx@>>>\kk$ are given by
$$\p'(Tx,y,X_1,\do,X_{r-1},Y_1,\do,Y_{r-1},R_1,\do,R_{r-1})=y|R_1,\do,R_{r-1}|,$$
$$\align&\ti h(Tx,y,X_1,\do,X_{r-1},Y_1,\do,Y_{r-1},R_1,\do,R_{r-1})\\&=
\sum_{j\in[1,r-1]}\la Y_j,R_j\ra
+\sum_{j\in[1,r-1]}\la{}_xA_j,u_j({}_yX_1,\do,{}_yX_j,Y_1,\do,Y_j,X_1,\do,X_j)\ra.\endalign$$

\subhead 1.5\endsubhead
Let $\tcx''$ be the variety of all
$$\align&(Tx,y,X_1,X_2,\do,X_{r-2},Y_1,Y_2,\do,Y_{r-2},R_1,R_2,\do,R_{r-1})\\&\in(T\bsl G)\T G\T
\fg^{(r-2)+(r-2)+(r-1)}\endalign$$
such that $xyx\i\in T$ and
$$u_j({}_yX_1,\do,{}_yX_j,Y_1,\do,Y_j,X_1,\do,X_j)\in{}_x\ft\text{ for }1\le j\le r'-1,$$
$$u_j({}_yX_1,\do,{}_yX_j,Y_1,\do,Y_j,X_1,\do,X_j)\in{}_x\fb\text{ if $j=r'$ and $r$ is odd }.$$
(The equations make sense since if $1\le j\le r'-1$ then $j\le r-2$ and since when $r$ is odd we have
$r'=r-r'-1\le r-2$.) 
We define $\mu:\tcx@>>>\tcx''$ by
$$\align&(Tx,y,X_1,X_2,\do,X_{r-1},Y_1,Y_2,\do,Y_{r-1},R_1,R_2,\do,R_{r-1})\m\\&
(Tx,y,X_1,X_2,\do,X_{r-2},Y_1,Y_2,\do,Y_{r-2},R_1,R_2,\do,R_{r-1}).\endalign$$
This is a vector bundle; for a fixed 
$$s=(Tx,y,X_1,X_2,\do,X_{r-2},Y_1,Y_2,\do,Y_{r-2},R_1,R_2,\do,R_{r-1})\in\tcx'',$$
 the fibre $\mu\i(s)$ can be identified with $\fg^2$ with coordinates $X_{r-1},Y_{r-1}$.
The restriction of $\ti h$ to $\mu\i(s)$ is of the form
$$(X_{r-1},Y_{r-1})\m\la Y_{r-1},R_{r-1}+{}_xA_{r-1}+c$$
where $c$ is a constant depending on $s$. We use that $\la{}_xA_j,{}_yX_{r-1}-X_{r-1}\ra=0$; this holds
since ${}^{yx\i}A_{r-1}={}^{x\i}A_{r-1}$ (recall that $xyx\i\in T$). Thus this restriction is affine linear
and is constant precisely when $R_{r-1}=-{}_xA_{r-1}$. Hence the results in 1.3 are applicable.
Let $\bar\cx$ be the variety of all
$$(Tx,y,X_1,X_2,\do,X_{r-2},Y_1,Y_2,\do,Y_{r-2},R_1,R_2,\do,R_{r-1})\in\tcx''$$
such that $R_{r-1}=-{}_xA_{r-1}$.
We define $\bar{\p}:\bar\cx@>>>G_r$, $\bar h:\bar\cx@>>>\kk$ by
$$\align&\bar{\p}(Tx,y,X_1,X_2,\do,X_{r-2},Y_1,Y_2,\do,Y_{r-2},R_1,R_2,\do,R_{r-1})\\&
=y|R_1,R_2,\do,R_{r-1}|,\endalign$$
$$\align&\bar h(Tx,y,X_1,X_2,\do,X_{r-2},Y_1,Y_2,\do,Y_{r-2},R_1,R_2,\do,R_{r-1})\\&=
\sum_{j\in[1,r-2]}\la Y_j,R_j\ra
+\sum_{j\in[1,r-2]}\la{}_xA_j,u_j({}_yX_1,\do,{}_yX_j,Y_1,\do,Y_j,X_1,\do,X_j)\ra\\&
+\la {}_xA_{r-1},u'_{r-1}({}_yX_1,\do,{}_yX_{r-2},Y_1,\do,Y_{r-2},X_1,\do,X_{r-2})\ra,\endalign$$
with notation of 1.1(a).
The inverse image of $\ce$ under $\bar\cx@>>>T$,
$$(Tx,y,X_1,X_2,\do,X_{r-2},Y_1,Y_2,\do,Y_{r-2},R_1,R_2,\do,R_{r-1})\m xyx\i$$
is denoted again by $\tce$. Then from 1.3(b) we deduce
$$\hK=\bar{\p}_!(\tce\ot\cl_{\bar h})[(r-5)\D].\tag a$$
Let $\cc\sub\fg$ be the $G$-orbit of $-A_{r-1}$ for the adjoint action, a regular semisimple orbit.
Let $V=\{(y,R)\in G\T\cc;{}^yR=R\}$. Let $\Si$ be the support of $\hK$ (a closed subset of $G_r$).
From (a) we see that
$$\Si\sub\{y|R_1,R_2,\do,R_{r-1}|\in G_r;(y,R_{r-1})\in V\}.$$
It is likely that $\Si$ is a smooth subvariety of $G_r$, isomorphic to a vector bundle over $V$ with fibres
isomorphic to $(\ft^\perp)^{r-2}$. We will we show that this is the case at least when $r\in\{2,3,4\}$. 
Moreover, it is likely that when $r$ is even, $\hK$ is up to shift the intersection cohomology complex 
associated to a local system of rank $1$ on the smooth closed subvariety $\Si$. We will show that this is the
case when $r\in\{2,4\}$ and that the analogous statement is not true when $r=3$.

\subhead 1.6\endsubhead
The method used in 1.5 to eliminate the variables $X_{r-1},Y_{r-1}$ can be used to eliminate all variables 
$X_{r-r'},\do,X_{r-1},Y_{r-r'},\do,Y_{r-1}$. 
Let $\tcx''_1$ be the variety of all
$$\align&(Tx,y,X_1,X_2,\do,X_{r-r'-1},Y_1,Y_2,\do,Y_{r-r'-1},R_1,R_2,\do,R_{r-1})\\&\in(T\bsl G)\T G\T
\fg^{r-1+2(r-r'-1)}\endalign$$
such that $xyx\i\in T$ and
$$u_j({}_yX_1,\do,{}_yX_j,Y_1,\do,Y_j,X_1,\do,X_j)\in{}_x\ft\text{ for }1\le j\le r'-1,$$
$$u_j({}_yX_1,\do,{}_yX_j,Y_1,\do,Y_j,X_1,\do,X_j)\in{}_x\fb\text{ if $j=r'$ and $r$ is odd }.$$
(The equations make sense since if $1\le j\le r'-1$ then $j\le r-r'-1$ and since when $r$ is odd we have
$r'=r-r'-1$.) 
We define $\mu_1:\tcx@>>>\tcx''_1$ by
$$\align&(Tx,y,X_1,X_2,\do,X_{r-1},Y_1,Y_2,\do,Y_{r-1},R_1,R_2,\do,R_{r-1})\m\\&
(Tx,y,X_1,X_2,\do,X_{r-r'-1},Y_1,Y_2,\do,Y_{r-r'-1},R_1,R_2,\do,R_{r-1}).\endalign$$
This is a vector bundle; for a fixed 
$$s=(Tx,y,X_1,X_2,\do,X_{r-r'-1},Y_1,Y_2,\do,Y_{r-r'-1},R_1,R_2,\do,R_{r-1})\in\bar{\cx},$$
the fibre $\mu\i(s)$ can be identified with $\fg^{2r'}$ with coordinates 
$$X_{r-r'},do,X_{r-1},Y_{r-r'},\do,Y_{r-1}.$$
The restriction of $\ti h$ to $\mu\i(s)$ in an affine linear function. This follows from the fact that
for $j\in[1,r-1]$, the Lie 
polynomial 
$$u_j(X'_1,\do,X'_j,Y_1,\do,Y_j,X_1,\do,X_j)$$
 is a linear combination of terms which are iterated 
brackets of indeterminates $X'_h,Y_h,X_h$ with sum of indices equal to $j$ (hence $\le r-1$) hence 
containing at most one $X'_h,Y_h$ or $X_h$ with $h\ge r-r'$. (If they contained more than one, we would have
$2(r-r')\le r-1$ hence $r\le 2r'-1$, a contradiction.) Hence the results in 1.3 are applicable and they 
result in a description of $\hK$ which does not involve $X_{r-r'},\do,X_{r-1},Y_{r-r'},\do,Y_{r-1}$. But
even after this method is applied, one needs further arguments to analyze $\hK$, as we will see in 
Sections 2 and 3.

\head 2. The cases $r=2$ and $r=4$\endhead
\subhead 2.1\endsubhead
In this subsection we assume that $r=2$. Now 
$$\bar\cx=\{(Tx,y,R_1)\in(T\bsl G)\T G\T\fg;xyx\i\in T, R_1=-{}_xA_1.\}$$
We have $\bar{\p}(Tx,y,R_1)=y|R_1|$ and $\bar h:\bar\cx@>>>\kk$ is identically $0$.
Using 1.5(a) we have $$\hK=\bar{\p}_!\tce[-3\D]\tag a$$
Note that $\bar{\p}$ defines an isomorphism of $\bar\cx$ with
$$\cz=\{y|R_1|\in G_2; R_1\in\cc,{}^yR_1=R_1\}$$
and that $\cz$ is closed in $G_2$ (we use that $\cc$ is closed in $\fg$).
Moreover $\cz$ is a smooth subvariety of $G_2$ and $\bar{\p}$ can be viewed as the imbedding $\cz@>>>G_2$.
Since $\cz$ is closed in $G_2$ and smooth, irreducible of dimension $\D$ we see that
$\bar{\p}_!\tce[\D]$ is a simple perverse sheaf on $G_2$. Hence $\hK[4\D]$
is a simple perverse sheaf on $G_2$  with support $\Si=\cz$. 
Using Laumon's theorem \cite{\LA}, it follows that 

(b) {\it $K[4\D]$ is a simple perverse sheaf on $G_2$.}

\subhead 2.2\endsubhead
We now assume (until the end of 2.5) that $r=4$. Now $\bar\cx$ is the variety of all
$$(Tx,y,X_1,X_2,Y_1,Y_2,R_1,R_2,R_3)\in(T\bsl G)\T G\T\fg^7$$
such that $xyx\i\in T$, ${}_yX_1-X_1+Y_1\in{}_x\ft$ and $R_3=-{}_xA_3$. We have 
$$\bar{\p}(Tx,y,X_1,X_2,Y_1,Y_2,R_1,R_2,R_3)=y|R_1,R_2,R_3|,$$
$$\align&\bar h(Tx,y,X_1,X_2,Y_1,Y_2,R_1,R_2,R_3)\\&
=\la Y_1,R_1\ra+\la Y_2,R_2\ra+\la {}_xA_1,{}_yX_1-X_1+Y_1\ra\\&+
\la {}_xA_2,{}_yX_2-X_2+Y_2+[{}_yX_1,Y_1]/2-[{}_yX_1,X_1]/2-[Y_1,X_1]/2\ra\\&
+\la {}_xA_3,[{}_yX_2,Y_1]+[X_2,X_1]-[{}_yX_2,X_1]-[Y_2,X_1]-[{}_yX_1,[{}_yX_1,Y_1]]/6\\&
-[Y_1,[{}_yX_1,Y_1]]/3+[X_1,[{}_yX_1,Y_1]]/2+[{}_yX_1,[{}_yX_1,X_1]]/6+[{}_yX_1,[Y_1,X_1]]/6 \\&
+[Y_1,[{}_yX_1,X_1]]/6+[Y_1,[Y_1,X_1]]/6-[X_1,[{}_yX_1,X_1]]/3-[X_1,[Y_1,X_1]]/3\ra.\endalign$$
We make a change of variable $Y_1=X_1-{}_yX_1+{}_x\t$ where $\t\in\ft$. Then $\bar\cx$
 becomes the variety of all
$$(Tx,y,\t,X_1,X_2,Y_2,R_1,R_2,R_3)\in(T\bsl G)\T G\T\ft\T\fg^6$$
such that $xyx\i\in T$ and and $R_3=-{}_xA_3$. Now $\bar{\p}:\bar\cx@>>>G_4$ and $\bar h:\bar\cx@>>>\kk$ become
$$\bar{\p}(Tx,y,\t,X_1,X_2,Y_2,R_1,R_2,R_3)=y|R_1,R_2,R_3|,$$
$$\align&\bar h(Tx,y,\t,X_1,X_2,Y_2,R_1,R_2,R_3)\\&
=\la X_1-{}_yX_1+{}_x\t,R_1\ra+\la Y_2,R_2\ra+\la {}_xA_1,{}_x\t\ra\\&+
\la{}_xA_2,{}_yX_2-X_2+Y_2+[{}_yX_1,X_1+{}_x\t]/2-[{}_yX_1,X_1]/2-[-{}_yX_1+{}_x\t,X_1]/2\ra\\&+
\la{}_xA_3,[{}_yX_2,X_1-{}_yX_1+{}_x\t]+[X_2,X_1]-[{}_yX_2,X_1]-[Y_2,X_1]\\&
-[{}_yX_1,[{}_yX_1,X_1+{}_x\t]]/6-[X_1-{}_yX_1+{}_x\t,[{}_yX_1,X_1+{}_x\t]]/3\\&
+[X_1,[{}_yX_1,X_1+{}_x\t]]/2+[{}_yX_1,[{}_yX_1,X_1]]/6+[{}_yX_1,[-{}_yX_1+{}_x\t,X_1]]/6\\&
 +[X_1-{}_yX_1+{}_x\t,[{}_yX_1,X_1]]/6+[X_1-{}_yX_1+{}_x\t,[-{}_yX_1+{}_x\t,X_1]]/6\\&
-[X_1,[{}_yX_1,X_1]]/3-[X_1,[-{}_yX_1+{}_x\t,X_1]]/3\ra.\endalign$$
For $i=1,2,3$ we have $[A_i,\t]=0$ since $\ft$ is abelian; it follows that $\la A_i,[\x,\t]\ra=0$ for any 
$\x\in\fg$. We also have $\la{}_xA_i,{}_yX_j-X_j\ra=0$; indeed the left hand side is
$\la {}^{yx\i}A_i-{}^{x\i}A_i,X_j\ra$ and this is zero since ${}^{xyx\i}A_i=A_i$ (recall that $xyx\i\in T$).
Similarly we have $\la{}_xA_3,[{}_yX_2,-{}_yX_1]+[X_2,X_1]\ra=0$; indeed, the left hand side is 
$\la{}^{yx\i}A_3-{}^{x\i}A_3,[X_1,X_2]\ra=0$. We see that
$$\align&\bar h(Tx,y,\t,X_1,X_2,Y_2,R_1,R_2,R_3)
=\la X_1-{}_yX_1+{}_x\t,R_1\ra\\&+\la Y_2,R_2\ra+
\la {}_xA_1,{}_x\t\ra+\la {}_xA_2,Y_2-[-{}_yX_1,X_1]/2\ra\\&
+\la {}_xA_3,-[Y_2,X_1]+[{}_yX_1,[{}_yX_1,{}_x\t]]/6+[X_1,[X_1,{}_x\t]]/6\\&
+[X_1,[{}_yX_1,{}_x\t]]/6+[{}_yX_1,[{}_yX_1,X_1]]/6+[X_1,[{}_yX_1,X_1]]/6\ra.\endalign$$
Next we use the identity 
$$\la{}_xA_3,[Z,[Z',{}_x\t]]\ra=\la {}_x\t,[Z,[Z',{}_xA_3]]\ra$$
for any $Z,Z'$ in $\fg$. (This follows from $[{}_xA_3,{}_x\t]=0$.) We also use the equality
$$\la{}_xA_3,[{}_yX_1,[{}_yX_1,{}_x\t]]\ra=\la{}_xA_3,[X_1,[X_1,{}_x\t]]\ra.$$
(Since ${}^{yx\i}\t={}^{x\i}\t$, ${}^{yx\i}A_3={}^{x\i}A_3$, the left hand side is
$$\la{}_xA_3,{}_y[X_1,[X_1,{}_x\t]]\ra=
\la{}^{yx\i}A_3,[X_1,[X_1,{}_x\t]]\ra=\la {}_xA_3,[X_1,[X_1,{}_x\t]]\ra,$$
as required.) We see that 
$$\align&\bar h(Tx,y,\t,X_1,X_2,Y_2,R_1,R_2,R_3)=\la Y_2,R_2+{}_xA_2-[X_1,{}_xA_3]\ra\\&
+\la{}_x\t,R_1+{}_xA_1+[X_1,[X_1,{}_xA_3]]/6+[X_1,[{}_yX_1,{}_xA_3]]/3 \ra+\la X_1-{}_yX_1,R_1\ra\\&
+\la{}_xA_2,[{}_yX_1,X_1]/2\ra+\la{}_xA_3,[{}_yX_1,[{}_yX_1,X_1]]/6+[X_1,[{}_yX_1,X_1]]/6\ra.\endalign$$

\subhead 2.3\endsubhead
Let $\ct=\{(Tx,y,X_1,R_1,R_2,R_3)\in(T\bsl G)\T G\T\fg^4;xyx\i\in T,R_3+{}_xA_3=0\}$.
Let $\ct_0$ be the closed subset of $\ct$ consisting of all $(Tx,y,X_1,R_1,R_2,R_3)$ such that 
$$R_2+{}_xA_2-[X_1,{}_xA_3]=0,$$
$$R_1+{}_xA_1+[X_1,[X_1,{}_xA_3]]/3+[X_1,[{}_yX_1,{}_xA_3]]/6\in({}_x\ft)^\perp.$$
Define 
$\ti h_0:\ct_0@>>>\kk$ by
$$\align&\ti h_0(Tx,y,X_1,R_1,R_2,R_3)=\la X_1-{}_yX_1,R_1\ra+\la {}_xA_2,[{}_yX_1,X_1]/2\ra\\&
+\la{}_xA_3,[{}_yX_1,[{}_yX_1,X_1]]/6+[X_1,[{}_yX_1,X_1]]/6\ra.\endalign$$
Define $\bar\cx@>\ph>>\ct@>\ph'>>G_4$ by
$$\ph(Tx,y,\t,X_1,X_2,Y_2,R_1,R_2,R_3)=(Tx,y,X_1,R_1,R_2,R_3),$$
$$\ph'(Tx,y,X_1,R_1,R_2,R_3)=y|R_1,R_2,R_3|$$
so that $\p'=\ph'\ph$. Now $\ph$ is a vector bundle with
fibres of dimension $N=2\D+\d$. Note that the restriction of $\ti h:\bar\cx@>>>\kk$ to any fibre of
$\ph$ is affine linear and this restriction is constant precisely at the fibres over points in $\ct_0$;
moreover the constant is given by the value of $\ti h_0$. Using 1.3(b), we see that
$$\hK=j_!(\tce\ot\cl_{\ti h_0})[-5\D-2\d]\tag a$$
where $j:\ct_0@>>>G_4$ is the restriction of $\ph'$.

\subhead 2.4\endsubhead
Let $R$ be a regular semisimple element in $\fg$.
Let $\ft_R$ be the centralizer of $R$ in $\fg$; let $T_R$ be the centralizer of $R$ in $T$.
For any $z\in T_R$ we define a linear map $\Xi_{R,z}:\ft_R^\perp@>>>\fg/\ft_R^\perp$ by
$\Xi_{R,z}(\x)=[X,{}^z\x]\mod\ft_R^\perp$ where
$X$ is any element of $\fg$ such that $\x=[X,R]$. Note that such $X$ exists; if $X'$ is a another element
such that $\x=[X',R]$, then ${}^{z\i}(X'-X)\in\ft_R$ hence $X'=X+\r$ for some $\r\in\ft_R$ and
$[X',{}^z\x]=[X,{}^z\x]+[\r,{}^z\x]$. Since $[\r,{}^z\x]\in\ft_R^\perp$ we see that our map $\Xi_{R,z}$ is
well defined.

\subhead 2.5\endsubhead
Let $\cc\sub\fg$ be the $G$-orbit of $-A_3$ for the adjoint action, a regular semisimple orbit.
Let $\cz$ be the subset of $G_4$ consisting of all $y|R_1,R_2,R_3|$ such that

$R_3\in\cc$;

$y\in T_{R_3}$;

$R_2+{}_xA_2\in\ft_{R_3}^\perp$ where $Tx\in T\bsl G$ is uniquely determined by $R_3=-{}_xA_3$;

$R_1+{}_xA_1+\Xi_{R_3,1}(R'_2)/3+\Xi_{R_3,y\i}(R'_2)/6=0\text{ in }\fg/\ft_{R_3}^\perp$ where 
$R'_2=R_2+{}_xA_2$.
\nl
Note that $\cz$ is closed in $G_4$ (we use that $\cc$ is closed in $\fg$). Moreover $\cz$ is a smooth 
subvariety of $G_4$. Indeed, $V=\{(y,R)\in G\T\cc;y\in T_R\}$ is clearly smooth and $\cz$ is a fibration 
over $V$ with fibres isomorphic to $\ft^\perp\T\ft^\perp$.

From the definitions we see that $\ct_0=\ph'{}\i\cz$ and that the restriction of $\ph'$ defines a morphism
$\un\ph':\ct_0@>>>\cz$ whose fibres are exactly the orbits of the free $\ft$-action on $\ct_0$ given by
$$\t:(Tx,y,X_1,R_1,R_2,R_3)\m(Tx,y,X_1+{}_x\t,R_1,R_2,R_3).$$
Clearly, the local system $\tce$ on $\ct_0$ is the inverse image under $\un\ph'$ of a local system on $\cz$
denoted again by $\tce$. Next we show that 

(a) {\it the function $\ti h_0:\ct_0@>>>\kk$ is constant on each orbit of the $\ft$-action on $\ct_0$}
\nl
that is, if $\t\in\ft$ and $(Tx,y,X_1,R_1,R_2,R_3)\in\ct_0$, then  
$$\ti h_0(Tx,y,X_1+{}_x\t,R_1,R_2,R_3)=\ti h_0(Tx,y,X_1,R_1,R_2,R_3).$$
Thus, we must show that
$$\align&\la X_1+{}_x\t-{}_yX_1-{}_x\t,R_1\ra+\la {}_xA_2,[{}_yX_1+{}_x\t,X_1+{}_x\t]/2\ra\\&
+\la{}_xA_3,[{}_yX_1+{}_x\t,[{}_yX_1+{}_x\t,X_1+{}_x\t]]/6+[X_1+{}_x\t,[{}_yX_1+{}_x\t,X_1+{}_x\t]]/6\ra\\&
=\la X_1-{}_yX_1,R_1\ra+\la{}_xA_2,[{}_yX_1,X_1]/2\ra\\&
+\la{}_xA_3,[{}_yX_1,[{}_yX_1,X_1]]/6+[X_1,[{}_yX_1,X_1]]/6\ra.\endalign$$
(We have used that ${}_{xy}\t={}_x\t$.) It is enough to show that
$$\align&\la{}_xA_2,[{}_x\t,X_1]/2\ra+\la{}_xA_2,[{}_yX_1,{}_x\t]/2\ra\\&
+\la{}_xA_3,[{}_yX_1,[{}_yX_1,{}_x\t]]/6+[{}_yX_1,[{}_x\t,X_1]]/6+[{}_x\t,[{}_yX_1,{}_x\t]]/6\\&
+[{}_x\t,[{}_yX_1,X_1]]/6+[{}_x\t,[{}_x\t,X_1]]/6+[X_1,[{}_yX_1,{}_x\t]]/6+[X_1,[{}_x\t,X_1]]/6\\&
+[{}_x\t,[{}_yX_1,X_1]]/6+[{}_x\t,[{}_yX_1,{}_x\t]]/6+[{}_x\t,[{}_x\t,X_1]]/6\ra=0.\endalign$$
Since $\la{}_xA_i,[{}_x\t,\x]\ra=0$ for any $\x\in\fg$, we see that it is enough to show
$$\align&\la{}_xA_3,[{}_yX_1,[{}_yX_1,{}_x\t]]/6+[{}_yX_1,[{}_x\t,X_1]]/6
+[X_1,[{}_yX_1,{}_x\t]]/6\\&+[X_1,[{}_x\t,X_1]]/6\ra=0.\endalign$$
It is enough to show the following two equalities:
$$\la{}_xA_3,[{}_yX_1,[{}_yX_1,{}_x\t]]+[X_1,[{}_x\t,X_1]]\ra=0,\tag b$$
$$\la{}_xA_3,[{}_yX_1,[{}_x\t,X_1]]+[X_1,[{}_yX_1,{}_x\t]]\ra=0.\tag c$$
The left hand side of (b) is
$$\la{}_xA_3,{}_y[X_1,[X_1,{}_x\t]]-[X_1,[X_1,{}_x\t]]\ra=
\la{}^{yx\i}A_3-{}^{x\i}A_3,[X_1,[X_1,{}_x\t]]\ra$$
and this is zero since ${}^{yx\i}A_3={}^{x\i}A_3$. The left hand side of (c) is
$$\la{}_xA_3,[{}_x\t,[{}_yX_1,X_1]]$$
and this is zero since $\la{}_xA_3,[{}_x\t,\x]\ra=0$ for any $\x\in\fg$. This proves (a). 

From (a) we see that there is a unique morphism
$\hat h:\cz@>>>\kk$ such that $\ti h_0(s)=\hat h(\un\ph'(s))$ for any $s\in\ct_0$.
It follows that $\cl_{\ti h_0}=\un\ph'{}^*\cl_{\hat h}$.
Now $j:\ct_0@>>>G_4$ in 2.3 is a composition $\un j\un\ph'$ where
$\un j:\cz@>>>G_4$ is the imbedding. It follows that 
$j_!(\tce\ot\cl_{\ti h_0})=\un j_!(\tce\ot\un\ph'_!\un\ph'{}^*\cl_{\hat h})\cong
\un j_!(\tce\ot \cl_{\hat h})[-2\d]$. Combining with 2.3(a) we see that
$$\hK\cong\un j_!(\tce\ot \cl_{\hat h})[-5\D-4\d].$$

Since $\cz$ is closed in $G_4$ and smooth, irreducible of dimension $3\D-2\d$ we see that
$\un j_!(\tce\ot\cl_{\hat h})[3\D-2\d]$ is a simple perverse sheaf on $G_4$. Hence $\hK[8\D+2\d]$ is a
simple perverse sheaf on $G_4$ with support $\Si=\cz$. Using Laumon's theorem \cite{\LA}, it follows that

(b) {\it $K[8\D+2\d]$ is a simple perverse sheaf on $G_4$.}

\head 3. The case $r=3$ \endhead
\subhead 3.1\endsubhead
In this section we assume that $r=3$. Now $\bar\cx$ is the variety of all
$$(Tx,y,X_1,Y_1,R_1,R_2)\in(T\bsl G)\T G\T\fg^4$$
such that $xyx\i\in T$, ${}_yX_1-X_1+Y_1\in{}_x\fb$ and $R_2=-{}_xA_2$. In our case we have
$$\hK=\bar{\p}_!(\tce\ot\cl_{\bar h})[-2\D]$$ 
where $\bar{\p}:\bar\cx@>>>G_3$ and $\bar h:\bar\cx@>>>\kk$ are given by
$$\bar{\p}(Tx,y,X_1,Y_1,R_1,R_2)=y|R_1,R_2|,$$
$$\align&\bar h(Tx,y,X_1,Y_1,R_1,R_2)=\la Y_1,R_1\ra+\la{}_xA_1,{}_yX_1-X_1+Y_1\ra\\&
+\la{}_xA_2,[{}_yX_1,Y_1]/2-[{}_yX_1,X_1]/2-[Y_1,X_1]/2\ra.\endalign$$
Let 
$$\align&\bar\cx'=\{(Tx,y,X_1,R_1,R_2,\b);(Tx,y)\in(T\bsl G)\T G,(X_1,R_1,R_2)\in\fg^3,\b\in{}_x\fb,\\&
xyx\i\in T,R_2=-{}_xA_2\}.\endalign$$
We define an isomorphism $\bar\cx@>\si>>\bar\cx'$ by 
$$(Tx,y,X_1,Y_1,R_1,R_2)\m(Tx,y,X_1,R_1,R_2,\b)$$
where $\b\in{}_x\fb$ is given by $\b={}_yX_1-X_1+Y_1$. We identify $\bar\cx=\bar\cx'$ via this isomorphism. 
Then $\bar{\p},\bar h$ become 
$$\bar{\p}(Tx,y,X_1,R_1,R_2,\b)=y|R_1,R_2|,$$ 
$$\align&\bar h(Tx,y,X_1,R_1,R_2,\b)=\la X_1-{}_yX_1+\b,R_1\ra+\la{}_xA_1,\b\ra\\&
+\la{}_xA_2,[{}_yX_1,X_1-{}_yX_1+\b]/2-[{}_yX_1,X_1]/2-[X_1-{}_yX_1+\b,X_1]/2\ra\\&
=\la X_1-{}_yX_1,R_1\ra+\la{}_xA_1+R_1+[{}_xA_2,X_1+{}_yX_1]/2,\b\ra+\la{}_xA_2,[{}_yX_1,X_1]/2\ra.
\endalign$$
Let 
$$Z=\{(Tx,y,X_1,R_1,R_2)\in(T\bsl G)\T G\T\fg^3;xyx\i\in T,R_2=-{}_xA_2\},$$
$$Z_0=\{(Tx,y,X_1,R_1,R_2)\in Z;{}_xA_1+R_1+[{}_xA_2,X_1+{}_yX_1]/2\in{}_x\fn\}$$
Define $\p'_0:Z_0@>>>G_3$, $\ti h_0:Z_0@>>>\kk$ by
$$\p'_0(Tx,y,X_1,R_1,R_2)=y|R_1,R_2|,$$ 
$$\ti h_0(Tx,y,X_1,R_1,R_2)=\la X_1-{}_yX_1,R_1\ra+\la{}_xA_2,[{}_yX_1,X_1]/2\ra.$$
The map $\bar\cx'@>>>Z$ given by $(Tx,y,X_1,R_1,R_2,\b)\m(Tx,y,X_1,R_1,R_2)$ is a vector bundle
with fibres isomorphic to $\fb$. Applying 1.3(b) to this vector bundle we see that
$$\hK=\p'_{0!}(\tce\ot\cl_{\ti h_0})[-3\D-\d]$$ 
where the inverse image of $\ce$ under $Z_0@>>>T$, $(Tx,y,X_1,R_1,R_2)\m xyx\i$ is denoted again by 
$\tce$. (We have used that $2\dim\fb=\D+\d$.) 

For any $R\in\cc$ (see 1.5) let $T_R$ be the centralizer of $R$ in $G$ and let $\ft_R$ be the centralizer of
$R$ in $\fg$. Let $\car\sub\Hom(\ft,\kk^*)$ be the set of roots of $\fg$ with respect to $\ft$; for any 
$\a\in\car$ let $\fg^\a$ be the corresponding ($1$-dimensional) root subspace and let $e^\a:T@>>>\kk^*$ be 
the correspondings root homomorphism.

Let $\car^+=\{\a\in\car;\fg^\a\sub\fn\}$, $\car^-=\car-\car^+$. For $R\in\cc$ let 
$\fg^-_R=\op_{a\in\car^-}{}_x\fg^\a$, $\fg^+-_R=\op_{a\in\car^+}{}_x\fg^\a$ (where $R=-{}_xA_2$); we have a 
direct sum decomposition $\fg=\fg^-_R\op\ft_R\op\fg^+_R$. Hence for any $X\in\fg$ we can write 
uniquely $X=X^-_R+X^0_R+X^+_R$ with $X^-_R\in\fg^-_R$, $X^0_R\in\ft_R$, $X^+_R\in\fg^+_R$. Let $\hZ$ be the 
variety of all $(y,X,R_1,R)$ where $R\in\cc$, $R_1\in\fg$, $y\in T_R$, $X\in\fg^-_R$ such that
${}_xA_1+R_1+[{}_xA_2,X+{}_yX]/2\in{}_x\fn$ for some/any $x\in G$ such that $R=-{}_xA_2$.
Define $\hat\p:\hZ@>>>G_3$, $\hat h:\hZ@>>>\kk$ by $\hat\p(y,X,R_1,R)=y|R_1,R|$,
$$\hat h(y,X,R_1,R)=\la X-{}_yX,R_1\ra.$$
We define $\z:Z_0@>>>\hZ$ by $(Tx,y,X,R_1,R)\m(y,X^-_R,R_1,R)$.
This is well defined since, if $\b\in\fb$, $R=-{}_xA_2$ and $y\in T_R$, then 
$[{}_xA_2,{}_x\b+{}_y({}_x\b)]\in{}_x\fn$.
Now $\z$ is a vector bundle with fibres isomorphic to $\fb$. Note also that $\p'_0=\hat\p\z$. We show that
$\ti h_0=\hat h\z$.

For a fixed $(Tx,y,X,R_1,R)\in Z_0$ we have ${}_xA_1+R_1+[{}_xA_2,X+{}_yX]/2\in{}_x\fn$ and in particular   
$$R_1^0+{}_xA_1=0,$$
$$R_1^-=-[{}_xA_2,X^-+{}_y(X^-)]/2 \tag a$$
where we write $X^+,X^-,X^0$ instead of $X^+_R,X^-_R,X^0_R$.
We must show that 
$$\la X-{}_yX,R_1\ra+\la{}_xA_2,[{}_yX,X]/2\ra=\la X^--{}_y(X^-),R_1\ra$$
or equivalently
$$\align&\la X^++X^0-{}_y(X^++X^0),R_1\ra+\la{}_xA_2,[{}_y(X^++X^0),X^-]/2\\&+[{}_y(X^-),X^++X^0]/2
+[{}_y(X^++X^0),X^++X^0]/2\ra=0,\endalign$$
that is,
$$\la X^+-{}_y(X^+),R_1^-\ra+\la{}_xA_2,[{}_y(X^+),X^-]/2+[{}_y(X^-),X^+]/2\ra=0.$$
In the left hand side we replace $R_1^-$ by the expression (a) and we obtain
$$\align&\la X^+-{}_y(X^+),-[{}_xA_2,X^-+{}_y(X^-)]/2\ra+\la{}_xA_2,[{}_y(X^+),X^-]/2+[{}_y(X^-),X^+]/2\ra\\&
=\la{}_xA_2,[X^-+{}_y(X^-),X^+-{}_y(X^+)]/2+[{}_y(X^+),X^-]/2+[{}_y(X^-),X^+]/2\ra\\&
=\la{}_xA_2,[X^-,X^+]-[{}_y(X^-),{}_y(X^+)]/2\ra
=\la {}^{yx\i}A_2-{}^{x\i}A_2,[X^-,X^+]\ra=0\endalign$$
since ${}^{yx\i}A_2-{}^{x\i}A_2=0$. Thus our claim is proved. Applying 1.3(b) to the vector bundle $\z$
we deduce 
$$\hK=K'[-4\D-2\d],\qua K'=\hat\p_!(\tce\ot\cl_{\hat h}),$$
where the inverse image of $\ce$ under $\hZ@>>>T$, $(y,X,R_1,R)\m xyx\i$ (where $R=-{}_xA_2$) 
is denoted again by $\tce$. (We have used that $2\dim\fb=\D+\d$.) 

\subhead 3.2\endsubhead
Let $\cz$ be the set of all $y|R_1,R|\in G_3$ such that $R\in\cc$, $y\in T_R$ and $(R_1)^0_R=-{}_xA_1$ 
where $R=-{}_xA_2$. This is clearly a closed, smooth subvariety of $G_3$; it is irreducible of dimension 
$2\D-\d$. For any $R\in\cc$ let $\cz_R$ be the inverse image of $R$ under the map $G_3@>>>\cc$, 
$y|R_1,R|\m R$.

Let $\ch^i_{y,R_1,R}$ be the stalk at $y|R_1,R|\in G_r$ of the $i$-th cohomology sheaf of $K'$, see 3.1.
We want to describe the vector spaces $\ch^i_{y,R_1,R}$. Note that $\ch^i_{y,R_1,R}=0$ unless 
$y|R_1R|\in\cz$; we now assume that this condition is satisfied. Using 
$G$-equivariance and the $G$-homogeneity of $\cc$, we see that we may also assume that $R=-A_2$ and we write
$\ch^i_{y,R_1}$ instead of $\ch^i_{y,R_1,R}$. We have 
$$\ch^i_{y,R_1}=H^i_c(\hat\p\i(y|R_1,R|),\tce\ot\cl_{\hat h}).$$
For any $X\in \fg$ we can write uniquely $X=X^0+\sum_{\a\in\car}X^\a$ where $X^0\in\ft$, $X^\a\in\fg^\a$.
Note that we have ${}_yX=X^0+\sum_\a e^\a(y\i)X^\a$.

Then $\hat\p\i(y|R_1R|)$ can be identified with the affine space
$$\{(X^{-\a})_{\a\in\car^+};\a(A_2)(1+e^\a(y))X^{-\a}/2=R_1^{-\a}\}.\tag a$$
The restriction of $\hat h$ to $\hat\p\i(y|R_1R|)$ becomes the affine linear function
$$(X^{-\a})_{\a\in\car^+}\m\sum_{a\in\car^+}(1-e^\a(y))\la X^{-\a},R_1^\a\ra.\tag b$$
We consider several cases.

(1) for some $\a\in\car^+$ we have $1+e^\a(y)=0$ and $R_1^{-\a}\ne0$;

(2) for any $\a\in\car^+$ such that $1+e^\a(y)=0$ we have $R_1^{-\a}=0$ but for some such $\a$ we have 
$R_1^\a\ne0$;

(3) for any $\a\in\car^+$ such that $1+e^\a(y)=0$ we have $R_1^{-\a}=0$ and $R_1^\a=0$;
\nl
In case (1), the affine space (a) is empty and $\ch^i_{y,R_1}=0$.

In case (2), the affine space (a) is nonempty and (b) is non-constant hence $\ch^i_{y,R_1}=0$.

For any $y\in T$ we set $\Xi_y=\{\a\in\car^+;1+e^\a(y)=0\}$. 
In case (3), the affine space (a) is nonempty of dimension equal to 
$\sha(\Xi_y)$ and (b) is constant, hence $\ch^i_{y,R_1}$ is $1$-dimensional if $i=2\sha(\Xi_y)$ and is $0$ 
$i\ne2\sha(\Xi_y)$.

\subhead 3.3\endsubhead
For any subset $\Xi$ of $\car^+$ let $T^\Xi=\{y\in T;\Xi_y=\Xi\}$ (the sets $T^\Xi$ form a partition 
of $T$). Note that $T^\emp$ is an open dense subset of $T$.
For $\Xi\sub\car^+$ let $\cz_R^\Xi$ be the set of all $y|R_1R|\in\cz_R$ such that 
$y\in T^\Xi$ and $R_1^\a=0$, $R_1^{-\a}=0$ for all $\a\in\Xi$.
The subsets $\cz_R^\Xi$ are clearly disjoint. Let $\cz'_R=\cz_R-\cup_{\Xi\sub\car^+}\cz_R^\Xi$. 
Note that for $y|R_1,R|\in\cz_R^\Xi$, 
$\ch^i_{y,R_1}$ is $1$-dimensional if $i=2\sha(\Xi)$ and is $0$ if $i\ne2\sha(\Xi)$.
Moreover, for $y|R_1,R|\in\cz'_R$, we have $\ch^i_{y,R_1}=0$ for all $i$.
We show that for any $\Xi\sub\car^+$ we have 

(a) {\it $\dim\cx_R^\Xi+2\sha(\Xi)\le\dim\cz_R$ with strict inequality unless $\Xi=\emp$.}
\nl
Indeed, we have $\dim\cx_R^\Xi=\dim T^\Xi+2\sha(\car^+-\Xi)$. On the other hand, 
$\dim\cz_R=\d+2\sha(\car^+)$. Thus (a) is equivalent to $\dim T^\Xi\le\d$, with strict inequality unless 
$\Xi=\emp$; this is obvious.

From (a) we see that $K'|_{\cz_R}$ satisfies half of the defining properties of an intersection cohomology
complex (the ones not involving Verdier duality). It follows that $K'|_{\cz}$ itself satisfies the same half
of the defining properties of an intersection cohomology complex; moreover $\Si$ (the support of $K'$) is 
equal to $\cz$.  Hence the perverse cohomology sheaves of $K'[2\D-\d]$ satisfy ${}^pH^i(K'[2\D-d])=0$ for 
$i>0$ and ${}^pH^0(K'[2\D-d])$ is a simple perverse sheaf on $G_3$. Since $\hK=K'[-4\D-2\d]$, it follows that
${}^pH^i(\hK[6\D+\d])=0$ for $i>0$ and ${}^pH^0(\hK[6\D+\d])$ is a simple perverse sheaf on $G_3$. 
Using Laumon's theorem \cite{\LA} we deduce:

(b) {\it ${}^pH^i(K[6\D+\d])=0$ for $i>0$ and ${}^pH^0(K[6\D+\d])$ is a simple perverse sheaf on $G_3$.}

\subhead 3.4\endsubhead
Let $\cz^\emp$ be the set of all $y|R_1,R|\in\cz$ such that for any $\a\in\car^+$ we have $e^\a(xyx\i)\ne-1$
(where $R=-{}_xA_2$); this is an open dense subset of $\cz$. We define $f:\cz^\emp@>>>\kk$ by 
$$\align&f(y|R_1,R|)=\sum_{\a\in\car^+}
\fra{2}{\a(A_2)}\fra{1-e^\a(xyx\i)}{1+e^\a(xyx\i)}\la({}_xR_1)^\a,({}_xR_1)^{-\a}\ra\\&
=\sum_{\a\in\car}\fra{2}{\a(A_2)}\fra{1}{1+e^\a(xyx\i)}\la({}_xR_1)^\a,({}_xR_1)^{-\a}\ra.\tag a\endalign$$
For $(y,X,R_1,R)\in\hat\p\i(\cz^\emp)$ we have
$$f(\hat\p(y,X,R_1,R))=\hat h(y,X,R_1,R).\tag b$$
To prove (b) we can assume that $R=-A_2$. We then have
$$\hat h(y,X,R_1,R)=\sum_{a\in\car^+}(1-e^\a(y))\la X^{-\a},R_1^\a\ra.$$
Replacing here $X^{-\a}$ by $\fra{2}{\a(A_2)(1+e^\a(y))}R_1^{-\a}$ we obtain
$$\hat h(y,X,R_1,R)=\sum_{a\in\car^+}(1-e^\a(y))\fra{2}{\a(A_2)(1+e^\a(y))}\la R_1^{-\a},R_1^\a\ra
=f(y|R_1,R|).$$
as required.
Since $\hat h$ is an isomorphism $\hat\p\i(\cz^\emp)@>\si>>\cz^\emp$ (by results in 3.3), we see that
$K'|_{\cz^\emp}$ is the rank $1$ local system $\tce\ot\cl_f$ on $\cz^\emp$ where the inverse image of $\ce$ 
under $\cz^\emp@>>>T$, $y|R_1,R|\m xyx\i$ (where $R=-{}_xA_2$) is denoted again by $\tce$.
It follows that the simple perverse sheaf ${}^pH^0(\hK[6\D+\d])$ on $G_3$ is associated to the local system 
$\tce\ot\cl_f$ on the locally closed smooth irreducible subvariety $\cz^\emp$ of $G_3$.

\subhead 3.5\endsubhead
It is likely that $K'[2\D-d]$ is a simple perverse sheaf on $G_3$. This would imply that
$K[6\D+\d]$ is a simple perverse sheaf on $G_3$. 

\head 4. A comparison of two complexes \endhead
\subhead 4.1\endsubhead
We preserve the assumptions in 1.4. 
Let $L$ be as in 0.2 where $f_i:\ft@>>>\kk$ is $\t\m\la A_i,\t\ra$ for $i=1,\do,r-1$. 
In this section we describe a strategy for showing that a shift of $L$ is isomorphic to $K$ in 1.4.

We define a sequence of algebraic varieties  $\cx_r,\cx_{r-1},\do,\cx_{r-2r'}$ as follows. For 
$i\in\{r-r',r-r'+1,\do,r\}$ let $\cx_i$ be the variety consisting of all

$(Tx,y,X_1,\do,X_{r-1},Y_1,\do,Y_{r-1})\in(T\bsl G)\T G\T\fg^{2r-2}$ 
\nl
such that $xyx\i\in B$ and

$u_j({}_yX_1,\do,{}_yX_j,Y_1,\do,Y_j,X_1,\do,X_j)\in{}_x\fb$ for $1\le j\le i-1$.  
\nl
For $i\in\{r-2r',r-2r'+1,\do,r-r'-1\}$ let $\cx_i$ be the variety of all

$(Tx,y,X_1,\do,X_{r-1},Y_1,\do,Y_{r-1})\in (T\bsl G)\T G\T\fg^{2r-2}$ 
\nl
such that $xyx\i\in T$ and

$u_j({}_yX_1,\do,{}_yX_j,Y_1,\do,Y_j,X_1,\do,X_j)\in{}_x\ft$ for $1\le j\le r-r'-i-1$,

$u_j({}_yX_1,\do,{}_yX_j,Y_1,\do,Y_j,X_1,\do,X_j)\in{}_x\fb$ for $r-r'-i\le j\le r-r'-1$.  
\nl
For $i=r,r-1,\do,r-2r'$ we have a diagram
$$G_r@<\p_i<<\cx_i@>h_i>>\kk$$
where $\p_i(Tx,y,X_1,\do,X_{r-1},Y_1,\do,Y_{r-1})=y|Y_1,Y_2,\do,Y_{r-1}|$,
$$\align&h_i(Tx,y,X_1,\do,X_{r-1},Y_1,\do,Y_{r-1})\\&
=\sum_{j\in[1,r-1]}\la{}_xA_j,u_j({}_yX_1,\do,{}_yX_j,Y_1,\do,Y_j,X_1,\do,X_j)\ra.\endalign$$  
We define $\io:\cx_i@>>>T$ by $\io(Tx,y,X_1,\do,X_{r-1},Y_1,\do,Y_{r-1})=d(xyx\i)$
where $d:B@>>>T$ is as in 0.1. Let $\tce=\io^*\ce$. The inverse image of $\tce$ under various maps to $\cx_i$
is denoted again by $\tce$.

Let $L_i=\p_{i!}(\tce\ot\cl_{h_i})\in\cd(G_r)$.
Let $H$ be the kernel of the obvious map $B_r@>>>T$, a connected unipotent group of dimension
$r(\D+\d)/2-\d$. Note that $\cx_r$ is a principal $H$-bundle over $\tG_r$ in 0.2. It follows that
$L_r\cong K[r(\D+\d)-2\d]$. On the other hand we have $L=L_{r-2r'}$. We would like to show that
$L\cong K[r(\D+\d)-2\d]$. It is enough to show that

(a) $L_r=L_{r-1}=\do=L_{r-2r'}$.
\nl
Note that $\cx_r\sub\cx_{r-r'}\sub\do\sub\cx_{r'}\supset\cx_{r'-1}\supset\do\supset \cx_{r-2r'}$.
For $r'\le i\le r-1$ let $\p'_i:\cx_i-\cx_{i+1}@>>>G_r$, $h'_i:\cx_i-\cx_{i+1}@>>>\kk$ be the 
restrictions of $\p_i,h_i$ to $\cx_i-\cx_{i+1}$; let $L'_i=p'_{i!}(\tce\ot\cl_{h'_i})\in\cd(G_r)$.

For $r-2r'+1\le i\le r'$ let $\p''_i:\cx_i-\cx_{i-1}@>>>G_r$, $h''_i:\cx_i-\cx_{i-1}@>>>\kk$
be the restrictions of $\p_i,h_i$ to $\cx_i-\cx_{i-1}$; let $L''_i=p''_{i!}(\tce\ot\cl_{h''_i})\in\cd(G_r)$.

From the definitions we have distinguished triangles $(L'_i,L_i,L_{i+1})$ (for $i=r',r'+1,\do,r-1$) and
$(L''_i,L_i,L_{i-1})$ (for $r-2r'+1\le i\le r'$). Hence (a) would follow from statements (b),(c) below.

(b) $L''_{r'}=L''_{r'-1}=\do=L''_{r-2r'+1}=0$.

(c) $L'_{r'}=L'_{r'+1}=\do=L'_{r-1}=0$.
\nl
Here is a strategy to prove (b),(c).

For $r-2r'+1\le i\le r'$ one should partition $\cx_i-\cx_{i-1}$ into pieces isomorphic to $\fg$
so that the restriction of $h''_i$ to each piece is a nonconstant affine linear function and the
restriction of $\tce$ is $\bbq$. 
For $r'\le i\le r-1$ one should partition $\cx_i-\cx_{i+1}$ into pieces isomorphic to an affine space
so that the restriction of $h'_i$ to each piece is a nonconstant affine linear function and the
restriction of $\tce$ is $\bbq$. 
This should give the desired result. In 4.2-4.5 we carry out this strategy in several cases which are
sufficient to deal with the cases where $r\in\{2,3,4\}$.

\subhead 4.2\endsubhead
In this subsection we show that

(a) $L''_{r'}=0$.
\nl
Note that $\cx_{r'}-\cx_{r'-1}$ is the set of all

$(Tx,y,X_1,\do,X_{r-1},Y_1,\do,Y_{r-1})\in(T\bsl G)\T G\T\fg^{2r-2}$ 
\nl
such that $xyx\i\in B-T$ and

$u_j({}_yX_1,\do,{}_yX_j,Y_1,\do,Y_j,X_1,\do,X_j)\in{}_x\fb$ for $1\le j\le r-r'-1$.  
\nl
Let $Z$ be the set of all

$(Tx,y,X_1,\do,X_{r-2},Y_1,\do,Y_{r-1})\in (T\bsl G)\T G\T\fg^{2r-3}$ 
\nl
such that $xyx\i\in B-T$ and

$u_j({}_yX_1,\do,{}_yX_j,Y_1,\do,Y_j,X_1,\do,X_j)\in{}_x\fb$ for $1\le j\le r-r'-1$.  
\nl
Now $\p''_{r'}$ is a composition $\cx_{r'}-\cx_{r'-1}@>a>>Z@>a'>>G_r$ where 

$a(Tx,y,X_1,\do,X_{r-1},Y_1,\do,Y_{r-1})=(Tx,y,X_1,\do,X_{r-2},Y_1,\do,Y_{r-1})$,

$a'(Tx,y,X_1,\do,X_{r-2},Y_1,\do,Y_{r-1})=y|Y_1,\do,Y_{r-1}|$.
\nl
It is enough to show that $a_!(\tce\ot\cl_{h''_{r'}})=0$. 
Clearly, $\tce$ is the inverse image under $a$ of a local system on $Z$ denoted again by $\tce$. Hence 
$$a_!(\tce\ot\cl_{h''_{r'}})=\tce\ot a_!(\cl_{h''_{r'}})$$
and it is enough to
show that $a_!(\cl_{h''_{r'}})=0$. It is also enough to show that for any 
$s=(Tx,y,X_1,\do,X_{r-2},Y_1,\do,Y_{r-1})\in Z$ we have $H^*_c(a\i(s),\cl_{h''_{r'}})=0$.
Now $a\i(s)$ may be identified with the affine space $\fg$ with coordinate $X_{r-1}$ and $h''_{r'}$ is of the
form $X_{r-1}\m\la{}_xA_{r-1},{}_yX_{r-1}-X_{r-1}\ra+c$ where $c$ is a constant (for fixed $s$). It is 
enough to show that the linear form
$$X_{r-1}\m\la{}_xA_{r-1},{}_yX_{r-1}-X_{r-1}\ra=\la{}^{yx\i}A_{r-1}-{}^{x\i}A_{r-1},X_{r-1}\ra$$
on $\fg$ is not identically zero. If it was identically zero, we would have 
${}^{yx\i}A_{r-1}={}^{x\i}A_{r-1}$ 
hence $xyx\i$ centralizes $A_{r-1}$ hence $xyx\i\in T$ contradicting $xyx\i\in B-T$. This proves (a).

\subhead 4.3\endsubhead
In this subsection we show that

(a) $L''_{r'-1}=0$ (assuming that $r=4$).
\nl
Note that $\cx_1-\cx_0$ is the set of all $(Tx,y,X_1,X_2,X_3,Y_1,Y_2,Y_3)\in(T\bsl G)\T G\T\fg^6$ 
such that $xyx\i\in T$ and ${}_yX_1-X_1+Y_1\in{}_x\fb-{}_x\ft$.

We have a free action of $\fg$ on $\cx_1-\cx_0$:
$$E:(Tx,y,X_1,X_2,X_3,Y_1,Y_2,Y_3)\m(Tx,y,X_1,X_2+E,X_3+[E,X_1],Y_1,Y_2,Y_3).$$
The orbit of $(Tx,y,X_1,0,X_3,Y_1,Y_2,Y_3)$ is 
$$\co=\{(Tx,y,X_1,E,X_3+[E,X_1],Y_1,Y_2,Y_3),E\in\fg\}.$$
It is enough to show that $H^*_c(\co,\tce\ot\cl_{h''_1})=0$ for any such $\co$. 
Clearly $\tce\cong\bbq$ on $\co$. Hence it is enough to show that $H^*_c(\co,\cl_{h''_1})=0$.
We can identify $\co$ 
with the affine space $\fg$ with coordinate $E$. On this affine space $h''_1$ is of the form
$$E\m\la{}_xA_2,{}_yE-E\ra+\la{}_xA_3,{}_y[E,X_1]-[E,X_1]+[{}_yE,Y_1]+[E,X_1]-[{}_yE,X_1]\ra+c$$
where $c$ is a constant (for our fixed $\co$). We have
$$\la{}_xA_2,{}_yE-E\ra=\la{}^{yx\i}A_2-{}^{x\i}A_2,E\ra=0$$
since ${}^{yx\i}A_2={}^{x\i}A_2$ (recall that $xyx\i\in T$). Hence $h''_1$ is of the form
$$E\m\la{}_xA_3,[{}_yE,{}_yX_1]+[{}_yE,Y_1]-[{}_yE,X_1]\ra+c=\la{}_yE,[\x,{}_xA_3]\ra+c$$
where $\x={}_yX-X_1+Y_1$. It is enough to show that the linear form
$E\m\la{}_yE,[\x,{}_xA_3]\ra$ on $\fg$ is not identically zero.
If it is identically zero we would have $[\x,{}_xA_3]=0$ that is, $\x$ is in the centralizer of
${}_xA_3$ so that $\x\in{}_x\ft$, contradicting $\x\in{}_x\fb-{}_x\ft$. This proves (a).

\subhead 4.4\endsubhead
In this subsection we show that

(a) $L'_{r-1}=0$.
\nl
Note that $\cx_{r-1}-\cx_r$ is the set of all

$(Tx,y,X_1,X_2,\do,X_{r-1},Y_1,Y_2,\do,Y_{r-1})\in (T\bsl G)\T G\T\fg^{2r-2}$ 
\nl
such that $xyx\i\in B$ and

$u_j({}_yX_1,\do,{}_yX_j,Y_1,\do,Y_j,X_1,\do,X_j)\in{}_x\fb$ for $j=1,2,\do,r-2$,

$u_j({}_yX_1,\do,{}_yX_j,Y_1,\do,Y_j,X_1,\do,X_j)\n{}_x\fb$ for $j=r-1$.
\nl
Now $\p'_{r-1}$ is a composition $\cx_{r-1}-\cx_r@>a>>Z@>a'>>G_r$ where $Z$ is the set of all
$$(Bx,y,X_1,X_2,\do,X_{r-1},Y_1,Y_2,\do,Y_{r-1})\in (B\bsl G)\T G\T\fg^{2r-2}$$
satisfying the same conditions as the points of $\cx_{r-1}-\cx_r$ and $a$ is the obvious map.
It is enough to show that $a_!(\tce\ot\cl_{h'_{r-1}})=0$. 
Clearly $\tce$ is the inverse image under $a$ of a local system on $Z$ 
denoted again by $\tce$. Hence $a_!(\tce\ot\cl_{h'_{r-1}})=\tce\ot a_!(\cl_{h'_{r-1}})$ and it is enough
to show that $a_!(\cl_{h'_{r-1}})=0$. It is also enough to show that for any 
$s=(Bx,y,X_1,X_2,\do,X_{r-2},Y_1,Y_2,\do,Y_{r-1})\in Z$ we have $H^*_c(a\i(s),\cl_{h'_{r-1}})=0$.
Now $a\i(s)$ may be identified with $U$ by 
$$u\m(Tux,y,X_1,X_2,\do,X_{r-1},Y_1,Y_2,\do,Y_{r-1})$$
where $x$ is a fixed representative of $Bx$. For $j\in[1,r-1]$ we set
$$\x_j={}^xu_j({}_yX_1,\do,{}_yX_j,Y_1,\do,Y_j,X_1,\do,X_j)\in\fg.$$
Then $h'_{r-1}$ becomes the function $U@>>>\kk$ given by
$$\align&u\m\sum_{j\in[1,r-1]}\la {}_{ux}A_j,{}_x\x_j\ra\\&
=\sum_{j\in[1,r-2]}(\la A_j,{}^u\x_j-\x_j\ra+\la A_j,\x_j\ra)
+\la A_{r-1},\x_{r-1}\ra+\la {}_uA_{r-1}-A_{r-1},\x_{r-1}\ra.\endalign$$
For $j\in[1,r-2]$ we have $\x_j\in\fb$ hence ${}^u\x_j-\x_j\in\fn$ so that
$\la A_j,{}^u\x_j-\x_j\ra=0$. Thus $h'_{r-1}$ becomes the function $U@>>>\kk$ given by
$$u\m\la{}_uA_{r-1}-A_{r-1},\x_{r-1}\ra+c$$
where $c$ is a constant (for fixed $s$). We identify $U$ with $\fn$ by $u\m{}_uA_{r-1}-A_{r-1}$. 
Then $h'_{r-1}$ becomes the function $\fn@>>>\kk$ given by $\z\m\la \z,\x_{r-1}\ra+c$. This function is
affine linear and nonconstant since $\x_{r-1}\n\fb=\fn^\perp$. 
It follows that $H^*_c(a\i(s),\cl_{h'_{r-1}})=0$ and (a) is proved.

\subhead 4.5\endsubhead
In this subsection we show that

(a) $L'_{2r'-2}=0$ (assuming that $r=4$).
\nl
Note that $\cx_2-\cx_3$ is the set of all

$(Tx,y,X_1,X_2,X_3,Y_1,Y_2,Y_3)\in(T\bsl G)\T G\T\fg^6$ 
\nl
such that $xyx\i\in B$, 

(b) ${}_yX_1-X_1+Y_1\in{}_x\fb$,

(c) ${}_yX_2-X_2+Y_2+[{}_yX_1,Y_1]/2-[{}_yX_1,X_1]/2-[Y_1,X_1]/2\n{}_x\fb$.
\nl
Let $Z=\{(Tx,y,Y_1,Y_2,Y_3)\in(T\bsl G)\T G\T\fg^3;xyx\i\in B\}$.
The inverse image of $\ce$ under $Z@>>>T$, $(Tx,y,Y_1,Y_2,Y_3)\m d(xyx\i)$ is denoted by $\tce_0$.

Now $\p'_2$ is a composition $\cx_2-\cx_3@>a>>Z@>a'>>G_r$ where 
$$\align&a(Tx,y,X_1,X_2,X_3,Y_1,Y_2,Y_3)=(Tx,y,Y_1,Y_2,Y_3),\\& a'(Tx,y,Y_1,Y_2,Y_3)=y|Y_1,Y_2,Y_3|
\endalign$$. 
We have $a^*\tce_0=\tce$. It is enough to show that $a_!(\tce\ot\cl_{h'_2})=0$ that is,
$\tce_0\ot a_!(\cl_{h'_2})=0$. Thus it is enough to prove that $a_!(\cl_{h'_2})=0$.
Hence it is enough to show that for any $s=(Tx,y,Y_1,Y_2,Y_3)\in(T\bsl G)\T G\T\fg^3$ 
we have $H^*_c(a\i(s),\cl_{h'_2})=0$. 

Let $\cg=\{|E,E',E''|\in G_4;E\in\fn,E'\in\fn,E''\in\fn\}$; this is a closed subgroup
of $G_4$.

We fix a representative $x$ in $Tx$ and we define a free $\cg$-action on $a\i(s)$ by
$$\align&|E,E',E''|:(Tx,y,X_1,X_2,X_3,Y_1,Y_2,Y_3)\m(Tx,y,X_1+{}_xE,\\&
X_2+{}_xE'+[{}_xE,X_1]/2,X_3+{}_xE''\\&
+[{}_xE',X_1] -[{}_xE,[{}_xE,X_1]]/6-[X_1,[{}_xE,X_1]]/3,Y_1,Y_2,Y_3).\endalign$$
We verify that this action is well defined (that is, the equations (b),(c) are preserved).
To show that (b) is preserved it is enough to verify that ${}_{xy}E-{}_xE\in{}_x\fb$ or that
${}^{xy\i x\i}E-E\in\fb$; this follows from $E\in\fn$, $xyx\i\in B$. 
To show that (c) is preserved it is enough to verify that
$$\align&{}_{xy}E'+[{}_{xy}E,{}_yX_1]/2-{}_xE'-[{}_xE,X_1]/2+[{}_{xy}E,Y_1]/2\\&
-[{}_yX_1,{}_xE]/2-[{}_{xy}E,X_1]/2-[{}_{xy}E,{}_xE]/2-[Y_1,{}_xE]/2\in{}_x\fb\endalign$$
(when (b) holds) or that
$$[{}_{xy}E,{}_yX_1-X_1+Y_1]/2+[{}_xE,{}_yX_1-X_1+Y_1]/2-[{}_{xy}E,{}_xE]/2
+{}_{xy}E'-{}_xE'\in{}_x\fb$$
and this follows from (b) and from ${}_{xy}E\in{}_x\fb$, ${}_xE\in{}_x\fb$,
${}_{xy}E'\in{}_x\fb$, ${}_xE'\in{}_x\fb$.

It is enough to show that for any $\cg$-orbit $\co$ in $a\i(s)$ we have $H^*_c(\co,\cl_{h'_2})=0$. 
We may identify $\co=\cg$ using a base point $(Tx,y,X_1,X_2,X_3,Y_1,Y_2,Y_3)\in\co$ (with a fixed
representative $x$ for $Tx$) and we identify $\cg=\fn^3$ using $[E,E',E'']\lra(E,E',E'')$.
Then $h'_2$ becomes a function $h'':\fn^3@>>>\kk$ of the following form
(we have substituted $Y_1=X_1-{}_yX_1+{}_x\b$ where $\b\in\fb$):
$$(E,E',E'')\m h''(E,E',E'')=\la{}_xA_1,\x_1\ra+\la{}_xA_2,\x_2+\x'_2\ra+\la{}_xA_3,\x_3+\x'_3+\x''_3\ra$$
where 
$$\x_1={}_x\b+{}_{xy}E-{}_xE,$$
$$\align&\x_2={}_yX_2+[{}_{xy}E,{}_yX_1]/2-X_2-[{}_xE,X_1]/2+Y_2\\&
+[{}_yX_1+{}_{xy}E,X_1-{}_yX_1+{}_x\b]/2-[{}_yX_1+{}_{xy}E,X_1+{}_xE]/2\\&
-[X_1-{}_yX_1+{}_x\b,X_1+{}_xE]/2,\endalign$$
$$\x'_2={}_{xy}E'-{}_xE'$$
$$\align&\x_3={}_yX_3-X_3+Y_3-[{}_{xy}E,[{}_{xy}E,{}_yX_1]]/6-[{}_yX_1,[{}_{xy}E,{}_yX_1]]/3\\&
+[{}_xE,[{}_xE,X_1]]/6+[X_1,[{}_xE,X_1]]/3\\&
+[{}_yX_2+[{}_{xy}E,{}_yX_1]/2,X_1-{}_yX_1+{}_x\b]+[X_2+[{}_xE,X_1]/2,X_1+{}_xE]\\&
-[{}_yX_2+[{}_{xy}E,{}_yX_1]/2,X_1+{}_xE]-[Y_2,X_1+{}_xE]\\&
-[{}_yX_1+{}_{xy}E,[{}_yX_1+{}_{xy}E,X_1-{}_yX_1+{}_x\b]]/6\\&
-[X_1-{}_yX_1+{}_x\b,[{}_yX_1+{}_{xy}E,X_1-{}_yX_1+{}_x\b]]/3\\&
+[X_1+{}_xE,[{}_yX_1+{}_{xy}E,X_1-{}_yX_1+{}_x\b]]/2\\&
+[{}_yX_1+{}_{xy}E,[{}_yX_1+{}_{xy}E,X_1+{}_xE]]/6\\&
+[{}_yX_1+{}_{xy}E,[X_1-{}_yX_1+{}_x\b,X_1+{}_xE]]/6 \\&
+[X_1-{}_yX_1+{}_x\b,[{}_yX_1+{}_{xy}E,X_1+{}_xE]]/6\\&
+[X_1-{}_yX_1+{}_x\b,[X_1-{}_yX_1+{}_x\b,X_1+{}_xE]]/6\\&
-[X_1+{}_xE,[{}_yX_1+{}_{xy}E,X_1+{}_xE]]/3\\&
-[X_1+{}_xE,[X_1-{}_yX_1+{}_x\b,X_1+{}_xE]]/3\ra,\endalign$$
$$\align&\x'_3={}_{xy}E''+[{}_{xy}E',{}_yX_1] -{}_xE''-[{}_xE',X_1]\\&
+[{}_{yx}E',X_1-{}_yX_1+{}_x\b]+[{}_xE',X_1]-[{}_{xy}E',X_1],\endalign$$
$$\x''_3=[{}_xE',{}_xE]-[{}_{xy}E',{}_xE].$$
It is enough to show that for any fixed $E',E''$ in $\fn$, the function
$E\m h''_1(E)=h''(E,E',E'')$ is affine linear and nonconstant.
Let
$$S={}_yX_2-X_2+Y_2+[{}_yX_1,Y_1]/2-[{}_yX_1,X_1]/2-[Y_1,X_1]/2.$$
A computation shows that 
$$\x_1-C_1\in{}_x\fn,\x_2-C_2\in{}_x\fn,\x_3-[{}_xE,S]-C_3\in{}_x\fn,\x'_3=C_4$$ 
where $C_1,C_2,C_3,C_4$ are vectors in $\fg$ independent of $E$. 
Moreover, $\x'_2\in{}_x\fn$, $\x''_3\in{}_x\fn$. Since $\la{}_xA_i,{}_x\fn\ra=0$, for some 
constant $c\in\kk$ we have
$$h''_1(E)=\la{}_xA_3,[{}_xE,S]\ra+c=\la S,[{}_xA_3,{}_xE]\ra+c.$$
In particular, $E\m h''_1(E)$ is affine linear on $\fn$. To show that it is nonconstant it is enough to 
show that $E\m\la S,[{}_xA_3,{}_xE]\ra$ is not identically zero. Assume that it is identically zero.
Since $E\m[A_3,E]$ is a vector space isomorphism $\fn@>\si>>\fn$ it would follow that 
$\la S,{}_x\ti E\ra=0$ for any $\ti E\in\fn$ hence $S\in{}_x(\fn^\perp)$ that is, $S\in{}_x\fb$. This 
contradicts the definition of $\cx_2-\cx_3$ and proves (a).

\subhead 4.6\endsubhead
In this subsection we assume that $r\in\{2,3,4\}$. From 4.2, 4.3, 4.4, 4.5 we see that 4.1(b),(c) hold. 
Hence 4.1(a) holds. Hence $L\cong K[2\D]$ if $r=2$, $L\cong K[3\D+\d]$ if $r=3$ and $L\cong K[4\D+2\d]$ if 
$r=4$.

Using now 2.1(b), 2.5(b), 3.3(b) we deduce the folowing result.

\proclaim{Theorem 4.7} (a) $L[r\D]$ is a simple perverse sheaf on $G_r$ provided that $r=2$ or $r=4$.

(b) If $r=3$ we have ${}^pH^i(L[r\D])=0$ for $i>0$ and ${}^pH^0(L[r\D])=0$ is a simple perverse sheaf on 
$G_r$.
\endproclaim
It is likely that in fact $L[r\D]$ is a simple perverse sheaf on $G_r$ for any $r\ge2$.
For $r=3$ this would follow if the truth of the statements in 3.5 could be established.

\enddocument